\documentclass[11pt]{amsart}
\usepackage{amsopn}
\usepackage{amsmath,amsthm,amssymb}

\newcommand{\nc}{\newcommand}

\nc{\n}{\noindent}
\nc{\Id}{\mathbf{1}}
\nc{\vs}{\vspace{8pt}}
\nc{\alt}{\raise1pt\hbox{$\bigwedge$}}
\nc{\ncp}{\nabla^\mathrm{CP}}
\nc{\nhc}{\nabla^\mathrm{HC}}
\nc{\wncp}{\widehat\nabla^\mathrm{CP}}
\nc{\ct}{\cos_\theta}
\nc{\st}{\sin_\theta}
\nc{\ctt}{\cos_{\theta/2}}
\nc{\stt}{\sin_{\theta/2}}
\nc{\ft}{\hbox{$\frac12$}}
\nc{\Aff}{\mathit{Aff}}

\nc{\vg}{\mathfrak{v} }
\nc{\wg}{\mathfrak{w} }
\nc{\zg}{\mathfrak{z} }
\nc{\ngo}{\mathfrak{n} }
\nc{\kg}{\mathfrak{k} }
\nc{\mg}{\mathfrak{m} }
\nc{\bg}{\mathfrak{b} }
\nc{\ggo}{\mathfrak{g} }
\nc{\ggob}{\overline{\mathfrak{g}} }
\nc{\sog}{\mathfrak{so} }
\nc{\sug}{\mathfrak{su} }
\nc{\spg}{\mathfrak{sp} }
\nc{\slg}{\mathfrak{sl} }
\nc{\glg}{\mathfrak{gl} }
\nc{\cg}{\mathfrak{c} }
\nc{\hg}{\mathfrak{h} }
\nc{\tg}{\mathfrak{t} }
\nc{\ug}{\mathfrak{u} }
\nc{\dg}{\mathfrak{d} }
\nc{\ag}{\mathfrak{a} }
\nc{\pg}{\mathfrak{p} }
\nc{\sg}{\mathfrak{s} }
\nc{\pca}{\mathcal{P}}
\nc{\nca}{\mathcal{N}}

\nc{\vp}{\varphi}
\nc{\ddt}{\frac{{\rm d}}{{\rm d}t}}

\nc{\SO}{{\tt SO}}
\nc{\Spe}{{\tt Sp}}
\nc{\SL}{{\tt SL}}
\nc{\SU}{{\tt SU}}
\nc{\Or}{{\tt O}}
\nc{\U}{{\tt U}}
\nc{\GL}{{\tt GL}}
\nc{\Se}{{\tt S}}
\nc{\CL}{{\tt CL}}
\nc{\Spin}{{\tt Spin}}
\nc{\Pin}{{\tt Pin}}

\nc{\RR}{{\mathbb R}}
\nc{\HH}{{\mathbb H}}
\nc{\CC}{{\mathbb C}}
\nc{\ZZ}{{\mathbb Z}}
\nc{\QQ}{{\mathbb Q}}
\nc{\FF}{{\mathbb F}}
\nc{\NN}{{\mathbb N}}
\nc{\GG}{{\mathbb G}}
\nc{\JJ}{{\mathbb J}}
\nc{\II}{{\mathbb I}}
\nc{\KK}{{\mathbb K}}
\nc{\DD}{{\mathbb D}}

\nc{\ad}{\operatorname{ad}}
\nc{\Ad}{\operatorname{Ad}}
\nc{\rank}{\operatorname{rank}}
\nc{\Irr}{\operatorname{Irr}}
\nc{\End}{\operatorname{End}}
\nc{\Aut}{\operatorname{Aut}}
\nc{\Inn}{\operatorname{Inn}}
\nc{\Der}{\operatorname{Der}}
\nc{\Ker}{\operatorname{Ker}}
\nc{\Iso}{\operatorname{I}}
\nc{\Le}{\operatorname{L}}
\nc{\tr}{\operatorname{tr}}
\nc{\dif}{\operatorname{d}\!}
\nc{\sen}{\operatorname{sen}}
\nc{\modu}{\operatorname{mod}}
\nc{\Ric}{\operatorname{R}}
\nc{\Sym}{\operatorname{Sym}}
\nc{\sca}{\operatorname{sc}}
\nc{\scalar}{{\sf s}}
\nc{\grad}{\operatorname{grad}}
\nc{\ricci}{\operatorname{r}}
\nc{\riccin}{\operatorname{Ric}}
\nc{\Lie}{\operatorname{L}}
\nc{\tang}{\operatorname{T}}
\nc{\tm}{\operatorname{TM}}
\nc{\aff}{\mathfrak{aff}}

\theoremstyle{plain}
\newtheorem{thm}{Theorem}[section]
\newtheorem{prop}[thm]{Proposition}

\theoremstyle{definition}

\theoremstyle{remark}
\newtheorem*{rem}{Remark}
\newtheorem*{rems}{Remarks}

\newcommand{\ri}{{\rm (i)}}
\newcommand{\rii}{{\rm (ii)}}
\newcommand{\riii}{{\rm (iii)}}
\newcommand{\riv}{{\rm (iv)}}

\title{Double products and hypersymplectic structures on $\RR^{4n}$ }
\author[A.~Andrada]{Adri\'an Andrada}
\author[I. Dotti]{Isabel G. Dotti}
\address{CIEM, FaMAF, Universidad Nacional de C\'ordoba, Ciudad Universitaria, (5000) C\'ordoba, Argentina}
\thanks{2000 {\it Mathematics Subject Classification.} Primary: 53C50, 53C15;
Secondary: 53C26, 53C30. \\
Both authors were partially supported
by CONICET, ANPCyT, SECyT-UNC and ACC (Argentina).}
\email{andrada@mate.uncor.edu}
\email{idotti@mate.uncor.edu}

\begin{document}

\begin{abstract}
In this paper we give a procedure to construct hypersymplectic structures on $\RR^{4n}$ beginning with affine-symplectic data on $\RR^{2n}$. These structures are shown to be invariant by a 3-step nilpotent double Lie group and the resulting metrics are complete and not necessarily flat. Explicit examples of this construction are exhibited.
\end{abstract}

\maketitle

\section{Introduction}

A hypersymplectic structure on a $4n$-dimensional manifold $M$ is given by $(J,E,g)$ where $J$, $E$ are endomorphisms of the tangent bundle of $M$ such that
\[ J^2=-\Id, \quad E^2=\Id, \quad JE=-EJ,\]
$g$ is a neutral metric (that is, of signature $(2n,2n)$) satisfying
\[g(X,Y)= g(JX,JY)=-g(EX,EY)\] for all $X,Y$ vector fields on $M$ and the associated $2$-forms
\[ \omega_1(X,Y)=g(JX,Y),\quad \omega_2(X,Y)=g(EX,Y),\quad \omega_3(X,Y)=g(JEX,Y)\]
are closed. Manifolds carrying a hypersymplectic structure have a rich geometry, the neutral metric is K\"ahler and Ricci flat and its holonomy group is contained in ${\rm Sp}(2n,\RR)$ (\cite{Hi}). Moreover, the Levi Civita connection is flat, when restricted to the leaves of the canonical foliations associated to the product structure given by $E$ (see \cite{A}). Metrics associated to a hypersymplectic structure are also called neutral hyperk\"ahler (see \cite{K}).

Hypersymplectic structures have significance in string theory. In \cite{OV}, $N=2$ superstring theory is considered, showing that the critical dimension of such a string is $4$ and that the bosonic part of the $N=2$ theory corresponds to self-dual metrics of signature $(2,2)$ (see also \cite{BGPPR} and \cite{Hu}).

The quotient construction proved to be a powerful method to construct symplectic and hyperk\"ahler structures on manifolds. According to \cite{Hi} this method cannot be always applied in the setting of hypersymplectic structures. Compact complex surfaces with neutral hyperkahler metrics are biholomorphic to either complex tori or primary Kodaira surfaces and both carry non flat neutral hyperkahler metrics, by results of Kamada (see \cite{K}). In higher dimensions,
hypersymplectic structures on a class of compact quotients of 2-step nilpotent Lie groups were exhibited in \cite{FPP} in their search of neutral Calabi-Yau metrics.

The purpose of this paper is to give a procedure to construct
hypersymplectic structures on $\RR^{4n}$ with complete and not
necessarily flat associated neutral metrics. The idea behind the
construction will be to consider the canonical flat
hypersymplectic structure on $\RR^{4n}$ and then translate it by
using an appropriate group acting simply and transitively on
$\RR^{4n}$.  This group will be a double Lie group
$(\RR^{4n},\,\RR^{2n}\times\{0\},\,\{0\}\times\RR^{2n})$
constructed from affine data on $\RR^{2n}$.  The most important
feature achieved by this procedure is that the associated neutral
metrics obtained will be complete and invariant by a 3-step
nilpotent group of isometries (we note that homogeneity does not
necessarily imply completeness in the pseudoriemannian setting.)
The degree of nilpotency will be related to the flatness of the
metric since we will show that the neutral metric is flat if and
only if the group is at most 2-step nilpotent.  Moreover, we
provide explicit examples of 3-step nilpotent Lie groups admitting
compact quotients and carrying invariant complete and non flat
hypersymplectic structures. The induced metric on the associated
nilmanifold will be neutral K\"ahler, complete, non-flat and Ricci
flat.

The outline of this paper is as follows. In \S 2 we give to
$\RR^{4n}$ a structure of a nilpotent Lie group. Starting with a
fixed symplectic structure $\omega$ on $\RR^{2n}$ which is
parallel with respect to a pair of  affine structures we form the
associated double Lie group
$(\RR^{4n},\,\RR^{2n}\times\{0\},\,\{0\}\times\RR^{2n})$  and show
that it is at most 3-step nilpotent. In \S 3 we consider canonical
symplectic structures on $\RR^{4n}$, constructed from the given
$\omega$ on $\RR^{2n}$ and show they are invariant by the group
constructed in \S 2. We analyze in \S 4 the geometry of the
homogeneous metric obtained by using the double Lie group
structure given to $\RR^{4n}$ to translate the standard inner
product of signature $(2n,2n)$ on $\RR^{2n}\oplus \RR^{2n}$. The
resulting metric is hypersymplectic (hence Ricci flat), complete
and not necessarily flat. Finally, in \S 5, we exhibit explicitly
flat and non flat complete neutral metrics on $\RR^{4n}$ which are
also K\"ahler and Ricci flat. Complete flat hypersymplectic
metrics are constructed on 2-step nilpotent Lie groups of
dimension $8n$ (\S 5.1) carrying also a closed special form in the
sense of (\cite{FPP}, section 2) and thus, the procedure developed
in \cite{FPP} may be applied to produce non flat neutral
Calabi-Yau metrics on the associated Kodaira manifolds. In \S 5.2,
complete non flat hypersymplectic metrics are exhibited on
$\RR^8$, where a particular example is given by
\begin{align*}
g = & \left(-\frac{3}{2}(x'_1+x'_3)^2+x'_2-x'_4\right)(\dif x_1+\dif x_3)^2 +2(x'_1+x'_3)(\dif x_1+\dif x_3)(\dif x_2-\dif x_4)-x'_1\dif x_1\dif x'_1\\
&  \quad -\dif x_1\dif x'_2+ \dif x_2\dif x'_1 + x'_3(\dif x_1\dif x'_3 +\dif x_3\dif x'_1)+ (x'_1 +2x'_3)\dif x_3\dif x'_3 -\dif x_3\dif x'_4+\dif x_4\dif x'_3
\end{align*}
with respect to coordinates $x_1,\ldots,x_4,x'_1,\ldots,x'_4$. Furthermore, metrics with similar properties can be obtained in higher dimensions.

\

\section{Group structure on $\RR^{4n}$ }

The main goal of this section will be to attach a 3-step nilpotent Lie group to  data $(\nabla, \nabla', \omega)$ where $\nabla,\;\nabla'$ are affine structures on $\RR^{2n}$ compatible with a symplectic structure $\omega$.
We shall begin by recalling some definitions which will be used throughout this
article.

An affine structure (or a left symmetric algebra structure) on
$\RR^{n}$ is given by a connection $\nabla$,  that is, a bilinear map $\nabla :\RR^{n}\times \RR^{n} \rightarrow
\RR^{n}$ satisfying the following conditions:
\begin{eqnarray}\label{flat} \nabla_xy= \nabla_yx,\\
\label{conmutan} \nabla_x\nabla_y = \nabla_y\nabla_x \end{eqnarray}
for all $x,y \in \RR^{n}$. If $\omega$
is a non-degenerate skew-symmetric bilinear form on $\RR^{n}$, the affine structure $\nabla$ is compatible with $\omega$ if
\begin{equation}\label{compatible} \omega(\nabla_xy, z)= \omega(\nabla_xz,y), \quad x,y,z\in \RR^n.\end{equation}
We notice that affine structures $\nabla$ on $\RR^{2n}$ compatible with $\omega$ satisfy a condition stronger than~(\ref{conmutan}), namely,
\begin{equation}\label{xy}
\nabla_x\nabla_y=0,
\quad x,y \in \RR^{2n}.\end{equation}
The last equation follows from
\begin{multline*}
\omega(\nabla_x\nabla_yz,w)= \omega(\nabla_wx,\nabla_zy)=-\omega(\nabla_z\nabla_wx,y)=\\
=-\omega(\nabla_w\nabla_zx,y)=
-\omega(\nabla_yw,\nabla_xz)=-\omega(\nabla_x\nabla_yz,w).\end{multline*}

Let $\nabla$ and $\nabla'$ be two affine connections on $\RR^{2n}$ compatible with $\omega$ and assume furthermore that $\nabla$ and $\nabla'$ satisfy the following compatibility condition
\begin{equation} \label{compatibility}
\nabla_x\nabla'_y=\nabla_y\nabla'_x
\end{equation}
for all $x,y\in\RR^{2n}$.

From (\ref{compatibility}) and the compatibility of the connections with $\omega$, we obtain the following:
\begin{equation}\label{compatibility'}
\nabla'_x\nabla_y=\nabla'_y\nabla_x
\end{equation}
and
\begin{equation}\label{commutativity}
\nabla_x\nabla'_{y}= -\nabla'_{y}\nabla_x.
\end{equation}
for all $x,y\in\RR^{2n}$. Indeed, (\ref{compatibility'}) follows from
\begin{multline*}
\omega(\nabla'_x\nabla_yz,w)= -\omega(\nabla_yz,\nabla'_xw)=-\omega(\nabla_y\nabla'_xw,z)=
\\=-\omega(\nabla_x\nabla'_yw,z)=-\omega(\nabla_xz,\nabla'_yw)= \omega(\nabla'_y\nabla_xz,w),\end{multline*}
and (\ref{commutativity}) follows from
\begin{multline*}
\omega(\nabla_x\nabla'_yz,w)= \omega(\nabla_y\nabla'_xz,w)= \omega(\nabla_wy,\nabla'_zx)
=-\omega(\nabla'_z\nabla_wy,x)= \\=-\omega(\nabla'_w\nabla_zy,x)=-\omega(\nabla'_xw,\nabla_yz)=-\omega(\nabla'_x\nabla_yz,w)=
-\omega(\nabla'_y\nabla_xz,w).\end{multline*}

\

We shall show in the next theorem that two affine structures $\nabla$ and $\nabla'$ on $\RR^{2n}$ satisfying (\ref{compatibility}) and (\ref{compatibility'}) give rise to a Lie group structure on the manifold $\RR^{4n}$ such that $(\RR^{4n},\,\RR^{2n}\times\{0\},\,\{0\}\times\RR^{2n})$ is a double Lie group. We recall that a {\em double Lie group} is given by a triple $(G,G_+,G_-)$ of Lie groups such that $G_+,\,G_-$ are Lie subgroups of $G$ and the product $G_+\times G_-\rightarrow G,\;(g_+,g_-)\rightarrow g_+g_-$ is a diffeomorphism (see \cite{LW}). The next result shows that the additional condition (\ref{compatible}) of $\nabla$ and $\nabla'$ with a fixed $\omega$ imposes restrictions on the Lie group obtained.

\begin{thm}\label{nilpotent}
Let $\nabla$ and $\nabla'$ be two affine structures on $\RR^{2n}$ compatible with a symplectic form $\omega$ and satisfying also (\ref{compatibility}). Then $\RR^{2n}\times\RR^{2n}$ with the product given by
\begin{equation}\label{product}
(x,x')\cdot(y,y')=(x+\alpha(x',y),\beta(x',y)+y')
\end{equation}
where $x,x',y,y'\in\RR^{2n}$ and
\begin{equation}\label{ab}
\alpha(x',y)=y+\nabla'_yx'-\frac{1}{2}\nabla'_y\nabla_yx',\quad \beta(x',y)=x'-\nabla_{x'}y-\frac{1}{2}\nabla_{x'}\nabla'_{x'}y
\end{equation}
is a 3-step nilpotent double Lie group. Furthermore, the associated Lie bracket on its Lie algebra $\RR^{2n}\oplus\RR^{2n}$ is
\begin{equation}\label{bracket} [(x,x'),(y,y')]= (\nabla'_yx'-\nabla'_xy', \nabla_xy'-\nabla_yx').
\end{equation}
\end{thm}

\begin{proof}
Let us set $\RR^{2n}_+:=\RR^{2n}\times\{0\}$ and $\RR^{2n}_-:=\{0\}\times\RR^{2n}$.
The maps $\alpha:\RR^{2n}_-\times\RR^{2n}_+\longrightarrow\RR^{2n}_+$ and $\beta:\RR^{2n}_-\times\RR^{2n}_+\longrightarrow\RR^{2n}_-$ satisfy the conditions
\begin{gather}
\alpha_0=\Id,\quad \alpha_{x'}(0)=0, \quad \alpha_{x'+y'}=\alpha_{x'}\circ\alpha_{y'}\\
\beta_0=\Id,\quad \beta_{y}(0)=0, \quad \beta_{x+y}=\beta_{y}\circ\beta_{x}
\end{gather}
for all $x',y'\in\RR^{2n}_-,\,x,y\in\RR^{2n}_+$, where we denote $\alpha_{x'}(y):=\alpha(x',y)$ and $\beta_y(x'):=\beta(x',y)$ for $x'\in\RR^{2n}_-,\,y\in\RR^{2n}_+$. The above relations show that $\alpha$ is a left action of $\RR^{2n}_-$ on $\RR^{2n}_+$ and $\beta$ is a right action of $\RR^{2n}_+$ on $\RR^{2n}_-$. These maps satisfy also the following compatibility conditions
\[ \alpha_{x'}(x+y)=\alpha_{x'}x+\alpha_{\beta_xx'}y,\qquad \beta_x(x'+y')=\beta_{\alpha_{y'}x}x'+\beta_xy'.\]
According to \cite{LW}, the product given in \eqref{product} defines a Lie group structure on $\RR^{4n}$ such that $(\RR^{4n},\RR^{2n}_+,\RR^{2n}_-)$ is a double Lie group. Note that the neutral element of this group structure is $(0,0)$ and the inverse of $(x,x')\in\RR^{4n}$ is $(\alpha(-x',-x),\beta(-x',-x))$.

Let us determine now the associated Lie algebra. Linearizing the above actions, we obtain representations $\mu:\RR^{2n}_-\longrightarrow \End(\RR^{2n}_+)$ and $\rho:\RR^{2n}_+\longrightarrow \End(\RR^{2n}_-)$ given by
\[ \mu_{x'}(y)=\left.\ddt\right|_0(\dif\alpha_{tx'})_0(y),\quad \rho_{y}(x')=\left.\ddt\right|_{0}(\dif\beta_{ty})_0(x').\]
For $\alpha$ and $\beta$ given in \eqref{ab}, we obtain that
\[ (\dif\alpha_{x'})_0=\Id+\nabla'_{x'},\quad (\dif\beta_{y})_0=\Id-\nabla_y.\]
Hence,
\[ \mu_{x'}(y)=\nabla'_{x'}y,\quad\rho_y(x')=-\nabla_y{x'},\]
showing that the bracket of its Lie algebra is the one given in \eqref{bracket}.

If $\ad_{(x,x')}, \, x,x'\in \RR^{2n}$ stands for the transformation given by (\ref{bracket}), then, using that $\nabla$ and $\nabla'$ are torsion-free (see (\ref{flat})) one has
\[  \ad_{(x,x')}= \begin{pmatrix}
     & & & & & \cr
     & \nabla'_{x'} & & & -\nabla'_x &\cr
     & & & & & \cr
     & & & & & \cr
     & -\nabla_{x'} & & & \nabla_x &\cr
     & & & & &\cr
\end{pmatrix} .  \]
From (\ref{xy}) applied to both $\nabla$ and $\nabla'$, we obtain
\[  \ad_{(x,x')}^2= \begin{pmatrix}
     & & & & & \cr
     & \nabla'_x \nabla_{x'} & & & -\nabla'_x \nabla_x&\cr
     & & & & & \cr
     & & & & & \cr
     & -\nabla_{x'}\nabla'_{x'} & & & \nabla_{x'} \nabla'_x&\cr
     & & & & &\cr
\end{pmatrix}.  \]
Finally, using (\ref{commutativity}),
\[ \ad_{(x,x')}^3=0.\]
Hence $\RR^{2n}\oplus\RR^{2n}$ is a 3-step nilpotent Lie algebra or $\RR^{2n}\times\RR^{2n}$ is a 3-step nilpotent Lie group, as claimed.
\end{proof}

\

We set the notation to be used in what follows. Since the construction of the Lie group structure on $\RR^{4n}$ in Theorem \ref{nilpotent} depends on the affine structures $\nabla$ and $\nabla'$, we will denote this Lie group by $N_{\nabla, \nabla'}$. The corresponding Lie algebra will be denoted $\ngo_{\nabla,\nabla'}$ and the abelian Lie subalgebras $\RR^{2n}\oplus\{0\},\,\{0\}\oplus\RR^{2n}$ will be denoted $\ngo_+,\,\ngo_-$, respectively. We note that $(\ngo_{\nabla,\nabla'},\ngo_+,\ngo_-)$ is a {\em double Lie algebra}, that is, $\ngo_+$ and $\ngo_-$ are Lie subalgebras of $\ngo_{\nabla,\nabla'}$ and $\ngo_{\nabla,\nabla'}=\ngo_+\oplus \ngo_-$ as vector spaces.

\

\begin{rems} $\ri$ The construction of the Lie algebra in Theorem \ref{nilpotent} can be made without requiring the affine structures to be compatible with a symplectic form. Indeed, if $\nabla$ and $\nabla'$ are affine structures on $\RR^m$ satisfying (\ref{compatibility}) and (\ref{compatibility'}), then the bracket (\ref{bracket}) defines a Lie algebra structure on $\RR^m\oplus\RR^m$ such that $\RR^{m}\oplus\{0\}$ and $\{0\}\oplus\RR^{m}$ are both abelian Lie subalgebras; hence $\RR^m\oplus\RR^m$ is 2-step solvable. Moreover, the centre of $\RR^m\oplus\RR^m$ is given by
\begin{equation}\label{centro}
\{(x,x') \in \RR^{m}\oplus \RR^{m}:\,\nabla_x=\nabla_x'=0,\;\nabla'_x=\nabla'_{x'}=0\}.
\end{equation}

\noindent $\rii$ If $\nabla$ is any affine structure on $\RR^m$ we denote by $A$ the associative (and commutative) algebra obtained from $\RR^m$ together with the product $x.y=\nabla_xy,\,x,y\in\RR^m$. In \cite{AS,BD}, $\aff(A)$ denoted the Lie algebra $A\oplus A$ with Lie bracket $[(a,b),(c,d)]=(0,ad-bc)$ with $a,b,c,d\in A$. Note that if, in $\ri$, we take $\nabla'=0$ then (\ref{compatibility}) and (\ref{compatibility'}) trivially hold, obtaining in this case a semidirect product which coincides with $\aff(A)$. This family of algebras and various geometric properties were considered in \cite{AS,BD}.
\end{rems}

\

\section{Invariant symplectic structures on $\RR^{4n}$}

In this section we show that $\ngo_{\nabla,\nabla'}$ carries three symplectic structures, obtained from the symplectic form $\omega$ in $\RR^{2n}$ compatible with $\nabla$ and $\nabla'$. These forms, defined at the Lie algebra level, give rise to left-invariant symplectic forms on the corresponding Lie group $N_{\nabla,\nabla'}$. Hence, $\RR^{4n}$ inherits symplectic structures which are invariant by this nilpotent group.

First, we recall that a symplectic structure on a Lie algebra $\ggo$ is a non-degenerate skew-symmetric bilinear form $\omega$ satisfying d$ \omega=0$, where
\begin{equation}\label{symplectic}  \text{d}\, \omega(x,y,z)= \omega (x,[y,z]) +\omega (y,[z,x]) +\omega (z,[x,y])
\end{equation}
for $x, y, z \in \ggo$.

A given symplectic form $\omega$ on $\RR^{2n}$ allows us to define the following non degenerate skew-symmetric bilinear forms on $\RR^{2n}\oplus\RR^{2n}$:
\begin{equation}\label{formas}
\begin{cases}
\omega_1((x,x'),(y,y')) & := \omega(x,y)+\omega(x',y'), \\
\omega_2((x,x'),(y,y')) & := -\omega(x,y')+\omega(y,x'),\\
\omega_3((x,x'),(y,y')) & := \omega(x,y)-\omega(x',y').
\end{cases}
\end{equation}
We show below that the above forms are closed with respect to the Lie bracket given in Theorem \ref{nilpotent}. Therefore, they define symplectic structures on $\ngo_{\nabla,\nabla'}$.

\begin{prop}
The 2-forms $\omega_1,\,\omega_2$ and $\omega_3$ are closed on $\ngo_{\nabla,\nabla'}$.
\end{prop}

\begin{proof}
Since $\RR^{2n}\oplus\{0\}$ and $\{0\}\oplus\RR^{2n}$ are abelian
subalgebras of $\ngo_{\nabla,\nabla'}$, it follows that the forms
$\omega_i,\,i=1,2,3,$ are closed if and only if
\[(\dif \omega_i)((x,0),(y,0),(0,z'))= (\dif
\omega_i)((0,x'),(0,y'),(z,0))=0 \] for all $x,y,z,x',y',z' \in
\RR^{2n}$. But
\begin{align*}
 (\dif \omega_i)((x,0),(y,0),(0,z')) & = \omega_i( [(y,0),(0,z')],(x,0))+ \omega_i( [(0,z'),(x,0)],(y,0))\\
& = \omega_i((-\nabla'_{y}z', \nabla_{y}z'), (x,0))+\omega_i((\nabla'_{x}z',-\nabla_{x}z'), (y,0))
\end{align*}
and
\begin{align*}
(\dif \omega_i)((0,x'),(0,y'),(z,0)) & = \omega_i( [(0,y'),(z,0)],(0,x'))+ \omega_i( [(z,0),(0,x')],(0,y'))\\
& =\omega_i((\nabla'_{z}y', -\nabla_{z}y'), (0,x'))+
\omega_i((-\nabla'_{z}x', \nabla_{z}x'), (0,y')).
\end{align*}
Using the expressions of $\omega_i,\,i=1,2,3$ given in (\ref{formas}), we compute
\[\begin{cases}
(\dif \omega_1)((x,0),(y,0),(0,z'))=(\dif \omega_3)((x,0),(y,0),(0,z'))=-\omega(-\nabla'_{y}z',x)+
\omega(\nabla'_{x}z',y), \\
(\dif \omega_1)((0,x'),(0,y'),(z,0))=-(\dif \omega_3)((0,x'),(0,y'),(z,0)) =
-\omega(\nabla_{z}y',x')+\omega(\nabla_{z}x',y') \end{cases}\]
and
\[ \begin{cases}
(\dif \omega_2)((x,0),(y,0),(0,z')) = \omega(x,\nabla'_{y}z')-\omega(y,\nabla_{x}z'),\\
(\dif \omega_2)((0,x'),(0,y'),(z,0)) = -\omega(\nabla'_{z}y',x')+
\omega(\nabla'_{z}x',y').\end{cases}\]
Since $\nabla$ and $\nabla'$ satisfy (\ref{flat}) and (\ref{compatible}), we obtain that $\dif\omega_i=0,\,i=1,2,3.$
\end{proof}

\

It follows from the definitions of the forms $\omega_i,\,i=1,2,3$ that:
\begin{enumerate}
\item[\ri] the restrictions of $\omega_1$ and $\omega_3$ to $\ngo_+$ and $\ngo_-$ are symplectic forms on these subalgebras;

\item[\rii] the Lie subalgebras $\ngo_+$ and $\ngo_-$ are Lagrangian subspaces of $\ngo_{\nabla,\nabla'}$ with respect to the symplectic form $\omega_2$.
\end{enumerate}

\

Let the form $\omega$ on $\RR^{2n}$ be given by $\omega=e^1\wedge e^2+e^3\wedge e^4+\cdots+e^{2n-1}\wedge e^{2n}$, where $\{e_1,\ldots,e_{2n}\}$ is a fixed basis of $\RR^{2n}$ and $\{e^1,\ldots,e^{2n}\}$ denotes the dual basis. Let us set $e_j:=(e_j,0)$ and $f_j:=(0,e_j),\,j=1,\ldots,2n$. Hence $\{e_1,\ldots,e_{2n},f_1,\ldots,f_{2n}\}$ is a basis of $\RR^{2n}\oplus\RR^{2n}$ and the forms $\omega_i,\,i=1,2,3$, can be written as
\begin{align*}
\omega_1 & = e^1\wedge e^2+\cdots+e^{2n-1}\wedge e^{2n}+f^1\wedge f^2+\cdots+f^{2n-1}\wedge f^{2n}, \\
\omega_2 & = -e^1\wedge f^2-\cdots-e^{2n-1}\wedge f^{2n}+e^2\wedge f^1+\cdots+e^{2n}\wedge f^{2n-1}, \\
\omega_3 & = e^1\wedge e^2+\cdots+e^{2n-1}\wedge e^{2n}-f^1\wedge f^2-\cdots-f^{2n-1}\wedge f^{2n}.
\end{align*}

\

\subsection{} Being $\ngo_{\nabla,\nabla'}$ a double Lie algebra, the endomorphism $E$ given by $E(x,y)=(x,-y)$ for $x,y\in\RR^{2n}$ is a {\em product structure} on $\ngo_{\nabla,\nabla'}$, that is, $E^2=\Id$ and $E$ is integrable, in the sense that it satisfies the condition
\begin{equation}\label{e}
E[(x,x'),(y,y')]=[E(x,x'),(y,y')]+[(x,x'),E(y,y')]-E[E(x,x'),E(y,y')].
\end{equation} for all $x,x',y,y'\in\RR^{2n}$. We note that the integrability of $E$ is equivalent to $\ngo_+$ and $\ngo_-$, the eigenspaces corresponding to the eigenvalues $\pm$ of $E$, being Lie subalgebras of $\ngo_{\nabla,\nabla'}$. Moreover, since $\ngo_+$ and $\ngo_-$ have equal dimension, $E$ is a {\em paracomplex structure} on $\ngo_{\nabla,\nabla'}$.

The symplectic form $\omega_2$ satisfies \[\omega_2(E(x,x'),E(y,y')) =-\omega_2((x,x'),(y,y'))\] for all $x,x',y,y'\in\RR^{2n}$. Therefore, $\{\ngo_{\nabla,\nabla'},\,E,\,\omega_2\}$ is an example of a {\em parak\"ahler Lie
algebra} in the sense of Kaneyuki (see \cite{Ka2}).

Another endomorphism on $\ngo_{\nabla,\nabla'}$ related to its decomposition as a double Lie algebra is given by $J(x,y)=(-y,x)$ for $x,y\in\RR^{2n}$. The endomorphism $J$ is a {\em complex structure} on $\ngo_{\nabla,\nabla'}$, that is, $J^2=-\Id$ and $J$ is integrable, i.e., it satisfies
\begin{equation}\label{j}
J[(x,x'),(y,y')]=[J(x,x'),(y,y')]+[(x,x'),J(y,y')]+J[J(x,x'),J(y,y')]
\end{equation}
for all $x,x',y,y'\in\RR^{2n}$. We note that $JE=-EJ$, and therefore $\{J,E\}$ is a complex product structure on $\ngo_{\nabla,\nabla'}$ (see \cite{AS}).

The symplectic form $\omega_1$ satisfies
\[ \omega_1(J(x,x'),J(y,y')) =\omega_1((x,x'),(y,y')) \]
for all $x,x',y,y'\in\RR^{2n}$. Hence, $\omega_1$ is a K\"ahler
form on $\ngo_{\nabla,\nabla'}$.

\

\section{Induced geometry on $\RR^{4n}$}\label{hypersymp}

In this section we analyze the properties of the metric on the manifold $\RR^{4n}$ obtained by left-translating by the Lie group $N_{\nabla, \nabla'}$, the standard inner product of signature $(2n,2n)$ on $\RR^{2n}\oplus\RR^{2n}$. We show that this metric is always complete and it is flat if and only if the Lie group $N_{\nabla, \nabla'}$ is 2-step nilpotent (see Theorem \ref{equiv} and Theorem \ref{completef}). Furthermore, this metric on $\RR^{4n}$ is hypersymplectic with respect to the structures $J$ and $E$ defined previously; in particular, it is neutral K\"ahler and Ricci-flat. Explicit examples will be given in subsequent sections.

Let us define a bilinear form $g$ on $\ngo_{\nabla,\nabla'}$ by
\begin{equation}
g((x,x'),(y,y'))= -\omega(x,y')+\omega(x',y)
\end{equation}
for all $(x,x'),(y,y')\in\ngo_{\nabla,\nabla'}$. It is clearly symmetric and non degenerate.
With respect to the basis $\{e_1,\ldots,e_{2n},f_1,\ldots,f_{2n}\}$, $g$ can be written as
\[ g=2\left(-e^1\cdot f^2-\cdots-e^{2n-1}\cdot f^{2n}+e^2\cdot f^1+\cdots+e^{2n}\cdot f^{2n-1}\right),\]
where $\cdot$ denotes the symmetric product of 1-forms. Moreover, $g$ satisfies the two following conditions:
\begin{gather}
g(J(x,x'),J(y,y'))=g((x,x'),(y,y'))\\
g(E(x,x'),E(y,y'))=-g((x,x'),(y,y'))
\end{gather}
for $x,x',y,y'\in\RR^{2n}$. Indeed,
\begin{align*}
g(J(x,x'),J(y,y')) & =g((-x',x),(-y',y))\\
                   & =-\omega(-x',y)+\omega(x,-y')\\
                   & =g((x,x'),(y,y'))
\end{align*}
and
\begin{align*}
g(E(x,x'),E(y,y')) & = g((x,-x'),(y,-y')) \\
                   & = \omega(x,y')-\omega(x',y) \\
                   & = -g((x,x'),(y,y')).
\end{align*}
Thus, $g$ is a Hermitian metric on $\ngo_{\nabla,\nabla'}$ with respect to both structures $J$ and $E$.

We note that the subalgebras $\ngo_+$ and $\ngo_-$ are both isotropic subspaces of $\ngo_{\nabla,\nabla'}$ with respect to $g$ and this metric has signature $(2n,2n)$.
Moreover, it is easy to verify that the $2$-forms $\omega_1,\omega_2$ and $\omega_3$ can be recovered from $g$ and the endomorphisms $J$ and $E$. Indeed we have
\begin{equation}\label{omega}
\begin{cases}
\omega_1((x,x'),(y,y')) =g(J(x,x'),(y,y')), \\
\omega_2((x,x'),(y,y')) =g(E(x,x'),(y,y')), \\
\omega_3((x,x'),(y,y')) =g(JE(x,x'),(y,y')).
\end{cases}
\end{equation}

\

The endomorphisms $J$ and $E$ of $\ngo_{\nabla,\nabla'}$, as well as the $2$-forms $\omega_1,\omega_2$ and $\omega_3$ and the metric $g$ can be extended to the group $N_{\nabla,\nabla'}$ by left translations. Hence, $N_{\nabla,\nabla'}$ is equipped with:
\begin{enumerate}
\item a complex structure $J$ and a product structure $E$ such that $JE=-EJ$;
\item a (pseudo) Riemannian metric $g$ such that
\begin{align*}
g(J(x,x'),J(y,y')) & = g((x,x'),(y,y')),\\
g(E(x,x'),E(y,y')) & = -g((x,x'),(y,y'))
\end{align*}
for all $x,x',y, y'\in\Gamma(\tang (N_{\nabla,\nabla'}))$;
\item three symplectic forms $\omega_1,\omega_2$ and $\omega_3$ which satisfy (\ref{omega}).
\end{enumerate}

To summarize, we have obtained

\begin{thm}
The nilpotent Lie group $N_{\nabla,\nabla'}$ carries a left-invariant hypersymplectic structure given by the $3$-tuple $\{J,E,g\}$.
\end{thm}

In particular, $(N_{\nabla,\nabla'},J,g)$ is a (neutral) K\"ahler manifold and $g$ is a Ricci-flat metric. Note also that $E$ is a product structure on $N_{\nabla,\nabla'}$ and hence, there is a decomposition of $\tang (N_{\nabla,\nabla'})$ into the Whitney sum of two involutive distributions of the same rank which are interchanged by $J$.

\begin{rem}
The leaves of both foliations given by $E$ are Lagrangian submanifolds of the symplectic manifold $(N_{\nabla,\nabla'},\omega_2)$. Hence, $N_{\nabla,\nabla'}$ is an example of a homogeneous parak\"ahler manifold (see \cite{Ka2}).
\end{rem}

\

\subsection{Curvature and completeness of $g$}

Since $g$ is left-invariant, the Levi-Civita connection $\nabla^g$ can be computed on left-invariant vector fields, i.e., on the Lie algebra $\ngo_{\nabla,\nabla'}$. After a computation one finds that
\begin{equation}\label{Levi-Civita}
{\nabla^g}_{(x,x')}= \begin{pmatrix}
    & & & & & \cr
    & \nabla_x +{\nabla'}_{x'} & & & 0 &\cr
    & & & & & \cr
    & & & & & \cr
    & 0 & & & \nabla_x + {\nabla'}_{x'} &\cr
    & & & & &\cr
\end{pmatrix} .  \end{equation}

One can verify, using the above expression of $\nabla^g$, that $J$ and $E$ are parallel with respect to the Levi-Civita connection.

We will show next that this connection need not be flat. If $R$ denotes the curvature of $\nabla^g$, it is easily seen (using (\ref{xy})) that $R((x,0),(y,0))= R((0,x'),(0,y'))=0$.  Moreover,
\[ R((x,0),(0,y'))= \nabla^g_{(x,0)}\nabla^g_{(0,y')}-\nabla^g_{(0,y')}\nabla^g_{(x,0)}
-\nabla^g_{(-\nabla'_xy',\nabla_xy')}\]
and using (\ref{Levi-Civita}) together with (\ref{flat}), (\ref{compatibility}) and (\ref{commutativity}) one  obtains
\[ R((x,0),(0,y'))=4\begin{pmatrix}
      & & & & & \cr
      & \nabla_x\nabla'_{y'} & & & 0 &\cr
      & & & & & \cr
      & & & & & \cr
      & 0 & & & \nabla_x \nabla'_{y'} & \cr
      & & & & &\cr
\end{pmatrix} = -4 \ad_{[(x,0),(0,y')]}.  \]
Since $R$ and the Lie bracket are skew-symmetric, one finally obtains
\[  R((x,x'), (y,y'))= -4\ad_{[(x,x')(y,y')]}, \]
thus showing that $\nabla^g$ will be flat if and only if $N_{\nabla,\nabla'}$ is $2$-step nilpotent. Note also that
\begin{equation}\label{R0} R=0 \quad \text{if and only if}\quad \nabla_x\nabla'_{y}=0\end{equation}
for all $x,y\in\RR^{2n}$. Summing up, we have shown

\begin{thm}\label{equiv}
The following conditions are equivalent:
\begin{enumerate}
\item[$\ri$] The Lie algebra $\ngo_{\nabla,\nabla'}$ is 2-step nilpotent;

\item[$\rii$] $\nabla_x\nabla'_{y}=0$ for all $x,y\in\RR^{2n}$;

\item[$\riii$] The hypersymplectic metric is flat.
\end{enumerate}
\end{thm}

\medskip

We end this section studying the completeness of $\nabla^g$. It follows from \cite{G} that $\nabla^g$ will be complete if and only if the differential equation on $\ngo_{\nabla,\nabla'}$
\begin{equation}\label{complete}
\dot{x}(t)=\ad^g_{x(t)}x(t)
\end{equation}
admits solutions $x(t)\in\ggo$ defined for all $t\in\RR$. Here $\ad_x^g$ means the adjoint of the transformation $\ad_x$ with respect to the metric $g$. It is easy to verify that the right hand side of (\ref{complete}) is given by $\ad^g_{x(t)}x(t)=-\nabla^g_{x(t)}x(t)$ for all $t$ in the domain of $x$ and thus we have to solve the equation
\begin{equation}\label{complete2}
\dot{x}(t)=-\nabla^g_{x(t)}x(t)
\end{equation}
The curve $x(t)$ on $\ngo_{\nabla,\nabla'}$ can be written as $x(t)=(a(t),b(t))$, where $a(t),b(t)\in\RR^{2n}$ are smooth curves on $\RR^{2n}$. Hence, using (\ref{Levi-Civita}), equation (\ref{complete2}) translates into the system
\begin{equation}\label{sistema} \begin{cases}
\dot{a}=-\nabla_aa-\nabla'_ba,\\
\dot{b}=-\nabla_ab-\nabla'_bb.
\end{cases} \end{equation}
Let us differentiate the first equation of the system above. We have
\begin{align*}
\ddot{a} & = -2\nabla_a\dot{a}-\nabla'_a\dot{b}-\nabla'_b\dot{a}\\
& = 2\nabla_a\nabla_aa+2\nabla_a\nabla'_ab+\nabla'_a\nabla_ab+\nabla'_a\nabla'_bb+\nabla'_b\nabla_aa+\nabla'_b\nabla'_ab\\
& = 0
\end{align*}
using (\ref{xy}), (\ref{compatibility}), (\ref{compatibility'}) and (\ref{commutativity}). In the same fashion, we differentiate the second equation of (\ref{sistema}) and obtain
\begin{align*}
\ddot{b} & = -\nabla_a\dot{b}-\nabla_b\dot{a}-2\nabla'_b\dot{b} \\
& = \nabla_a\nabla_ab+\nabla_a\nabla'_bb+\nabla_b\nabla_aa+\nabla_b\nabla'_ab+2\nabla'_b\nabla_ab+2\nabla'_b\nabla'_bb\\
& = 0
\end{align*}
using again (\ref{xy}), (\ref{compatibility}), (\ref{compatibility'}) and (\ref{commutativity}). Thus, there exist constant vectors $A,B,C,D\in\RR^{2n}$ such that
\[ a(t)=At+B,\qquad b(t)=Ct+D. \]
The explicit solution of the system (\ref{sistema}) with initial condition $x(0)=(a_0,b_0)$ is given by
\[ a(t)=(-\nabla_{a_0}a_0-\nabla'_{a_0}b_0)t+a_0, \qquad b(t)=(-\nabla_{a_0}b_0-\nabla'_{b_0}b_0)t+b_0. \]
Therefore, $x(t)$ is defined for all $t\in\RR$ and, in consequence, $\nabla^g$ is complete. Thus, we have obtained

\medskip

\begin{thm}\label{completef}
Hypersymplectic metrics on $\ngo_{\nabla,\nabla'}$ are always complete.
\end{thm}

\medskip

\begin{rem}
The completeness of $\nabla^g$ could have been dealt with in the following manner when $\nabla^g$ is flat. In this case, from results in \cite{S} we obtain that the completeness follows if the transformation $(x,x')\rightarrow {\nabla^g}_{(x,x')}(y,y')$ is nilpotent for every $(y,y')$. But, using (\ref{Levi-Civita}) one finds
\begin{equation}
{\nabla^g}{(y,y')}= \begin{pmatrix}
  & & & & & \cr
  & \nabla_y & & & {\nabla'}_{y} & \cr
  & & & & & \cr
  & & & & & \cr
  & \nabla_{y'} & & & {\nabla'}_{y'} & \cr
  & & & & &\cr
\end{pmatrix}. \end{equation}
Using (\ref{xy}) and (\ref{R0}) one shows that ${\nabla^g}{(y,y')}^2=0$ and hence the completeness follows.
\end{rem}

\

\section{Explicit examples}

In this section we will consider particular cases of the
constructions given previously.

In the first one we give explicit affine structures $\nabla$ and
$\nabla'$ such that the resulting group is $8n$-dimensional,
2-step nilpotent and with a $4n$-dimensional centre invariant by
the complex structure. It also admits compact quotients, hence the
associated nilmanifolds are Kodaira manifolds (see \cite{FPP}). In
this case the hypersymplectic metric obtained on the group is
complete and flat, hence it is isometric to the standard one. On
the other hand this Lie group carries a closed special form in the
sense of (\cite{FPP}, section 2) and thus, the procedure developed
in \cite{FPP} may be applied to produce non flat K\"ahler
Ricci-flat metrics on this family of Kodaira manifolds.

In the second one we give a 3-parameter family $\nabla$ and
${\nabla'}_{a,b,c}$ on $\RR^{4}$ satisfying all the requirements
to obtain a 3-step nilpotent group structure on $\RR^{8}$. In this
case the hypersymplectic metric obtained on $\RR^{8}$ will be
complete and non flat.  Moreover, the bracket relations will show
that compact quotients can be obtained. This example can be
generalized to higher dimensions, thus obtaining complete non flat
hypersymplectic metrics on $\RR^{4n},\,n\geq 2$, invariant by a
3-step nilpotent Lie group.

\subsection{Neutral K\"ahler Einstein metrics on Kodaira manifolds}

Let us consider on $\RR^{4n}=\text{span}\{e_1,\ldots,e_{4n}\}$ the
flat torsion-free connections $\nabla$ and $\nabla'$ given by
\[ \begin{cases}
      \nabla_{e_i}e_i=e_{i+1},\quad i \text{ odd, }1\leq i \leq 2n, \\
      \nabla_{e_j}e_k=0, \qquad \text{otherwise}
   \end{cases} \]
and
\[ \begin{cases}
      \nabla'_{e_i}e_i=e_{i+1},\quad i \text{ odd, }2n+1\leq i \leq 4n, \\
      \nabla'_{e_j}e_k=0, \qquad \text{otherwise}.
   \end{cases} \]
Clearly, $\nabla\nabla'=0=\nabla'\nabla$. Also, it is easy to see
that both connections are compatible with the standard symplectic
form $\omega$ on $\RR^{4n}$ given by $\omega=e^1\wedge
e^2+e^3\wedge e^4+\cdots+e^{4n-1}\wedge e^{4n}$.

We can form then the $8n$-dimensional Lie algebra
$\ngo_{\nabla,\nabla'}$ as in previous sections. It has a basis
$\{e_i,f_i:i=1,\ldots,4n\}$ where $e_i=(e_i,0)$ and $f_i=(0,e_i)$
and the only non-zero Lie brackets are
\[ [e_i,f_i]= \begin{cases}
                 f_{i+1}, \qquad i \text{ odd, }1\leq i \leq 2n, \\
                -e_{i+1},\, \quad i \text{ odd, }2n+1\leq i \leq 4n.
              \end{cases} \]
Let $N_{\nabla,\nabla'}$ denote the simply connected Lie group
associated to the Lie algebra $\ngo_{\nabla,\nabla'}$. Since the
structure constants are $0,\,1$ or $-1$, by Malcev's theorem
\cite{M}, there exists a discrete subgroup $\Gamma$ of
$N_{\nabla,\nabla'}$ such that
$M_{\Gamma}:=N_{\nabla,\nabla'}/\Gamma$ is compact. Using
(\ref{centro}), we obtain that the centre $\zg$ of
$\ngo_{\nabla,\nabla'}$ is given by
$\zg=\text{span}\{e_i,f_i:\,i\text{ is even}\}$, showing in
particular that this Lie group is 2-step nilpotent and $\zg$ is
$4n$-dimensional and stable under the action of $J$. Therefore,
the nilmanifold $M_{\Gamma}$ is an $8n$-dimensional Kodaira
manifold (see \cite{FPP}).

According to the results in \S\ref{hypersymp},
$N_{\nabla,\nabla'}$ carries a left-invariant hypersymplectic
structure, which is flat since $N_{\nabla,\nabla'}$ is $2$-step
nilpotent (see Theorem \ref{equiv}). Besides, since the centre
$\zg$ is, with respect to $\omega_1$, a Lagrangian subspace, the
symplectic form $\omega_1$ on the Lie group induces a closed
special 2-form on $M_{\Gamma}$. The method described in \cite{FPP}
can be applied in this case to produce Ricci-flat neutral K\"ahler
metrics on $M_\Gamma$.

We note that $\zg$ is a special Lagrangian subspace of $\ngo_{\nabla,\nabla'}$ with respect to the $J$-holomorphic form $\Phi=(\omega_2+i\omega_3)^{2n}$. This gives rise to special Lagrangian submanifolds on the symplectic manifold $M_\Gamma$.

The Levi-Civita connection $\nabla^g$ of the hypersymplectic metric $g$ on the group $N_{\nabla,\nabla'}$ is given by
\[ \begin{cases}
     \nabla^g_{e_i}e_i=e_{i+1}, \quad i \text{ odd, }1\leq i \leq 2n, \\
     \nabla^g_{e_i}f_i=f_{i+1}, \quad i \text{ odd, }1\leq i \leq 2n, \\
     \nabla^g_{f_i}e_i=e_{i+1}, \quad i \text{ odd, }2n+1\leq i \leq 4n, \\
     \nabla^g_{f_i}f_i=f_{i+1}, \quad i \text{ odd, }2n+1\leq i \leq 4n,
   \end{cases} \]
and $0$ in all the other possibilities. Note that we have the relations $\zg=\text{span}\{\nabla^g_xy:x,y\in\ngo_{\nabla,\nabla'}\}$ and $\nabla^gx=0$ for $x\in\zg$. Hence, Proposition 4.1 in \cite{FPP} can be applied, obtaining neutral Calabi-Yau metrics on compact quotients of the cotangent bundle of $N_{\nabla,\nabla'}$.

\

\subsection{Complete, non flat, neutral K\"ahler Einstein metrics on $\RR^{8}$} We will consider next $\RR^{4}=\text{span}\{e_1,\ldots,e_{4}\}$ equipped with two affine structures $\nabla$ and $\nabla'$ given by
\[ \nabla_{e_1}=\nabla_{e_3}=\begin{pmatrix}
                                      1 & 0 & 1 & 0 \cr
                                      0 & -1 & 0 & 1 \cr
                                     -1 & 0 & -1 & 0 \cr
                                      0 & -1 & 0 & 1 \cr
                             \end{pmatrix}     , \]
\[ \nabla_{e_2}=\begin{pmatrix}
                          0 & 0 & 0 & 0 \cr
                          -1 & 0 & -1 & 0 \cr
                          0 & 0 & 0 & 0 \cr
                          -1 & 0 & -1 & 0 \cr
                       \end{pmatrix} ,
\quad
\nabla_{e_4}=\begin{pmatrix}
                          0 & 0 & 0 & 0 \cr
                          1 & 0 & 1 & 0 \cr
                          0 & 0 & 0 & 0 \cr
                          1 & 0 & 1 & 0 \cr
                       \end{pmatrix}, \]
and
\[ \nabla'_{e_1}= \begin{pmatrix}
                          a & 0 & a & 0 \cr
                          b & -a & c & a \cr
                          -a & 0 & -a & 0 \cr
                          c & -a & -b+2c & a \cr
                       \end{pmatrix},
\quad
\nabla'_{e_2}= \begin{pmatrix}
                          0 & 0 & 0 & 0 \cr
                          -a & 0 & -a & 0 \cr
                          0 & 0 & 0 & 0 \cr
                          -a & 0 & -a & 0 \cr
                       \end{pmatrix}, \]
\[ \nabla'_{e_3}= \begin{pmatrix}
                          a & 0 & a & 0 \cr
                          c & -a & -b+2c & a \cr
                          -a & 0 & -a & 0 \cr
                          -b+2c & -a & -2b+3c & a \cr
                       \end{pmatrix},
\quad
\nabla'_{e_4}= \begin{pmatrix}
                          0 & 0 & 0 & 0 \cr
                          a & 0 & a & 0 \cr
                          0 & 0 & 0 & 0 \cr
                          a & 0 & a & 0 \cr
                       \end{pmatrix} \]
for $a,b,c\in\RR$. One can verify easily that $\nabla$ and $\nabla'$ satisfy equation (\ref{compatible}) with respect to the symplectic form $\omega=e^1\wedge e^2+e^3\wedge e^4$ on $\RR^{4}$. In order to see that the compatibility condition (\ref{compatibility}) holds, we observe that
\begin{equation}\label{ultima}
\nabla_{e_j}\nabla'_{e_k}= \begin{pmatrix}
                               0 & 0 & 0 & 0 \cr
                            -b+c & 0 & -b+c & 0 \cr
                               0 & 0 & 0 & 0 \cr
                            -b+c & 0 & -b+c & 0 \cr
                              \end{pmatrix}\end{equation}
for $(j,k)=(1,1),(1,3),(3,1),(3,3)$ and $\nabla_{e_j}\nabla'_{e_k}=0$ otherwise. From results in \S2, we may construct the 8-dimensional nilpotent Lie group $N_{\nabla,\nabla'}$ (see Theorem \ref{nilpotent}), whose underlying manifold is $\RR^4\times\RR^4$. The group structure is as follows
\[ (x,x')\cdot(y,y')=(x+\alpha(x',y),\beta(x',y)+y'),\]
where $\alpha$ and $\beta$ are given, in terms of their components, by
\begin{eqnarray*}
\alpha_1(x',y) & = & y_1+a(y_1+y_3)(x'_1+x'_3)\\
\alpha_2(x',y) & = & y_2+(by_1-ay_2+cy_3+ay_4)x'_1-a(y_1+y_3)x'_2\\
& & \quad +\left(cy_1-ay_2+(-b+2c)y_3+ay_4\right)x'_3 \\
& & \quad +a(y_1+y_3)x'_4+\ft(-b+c)(y_1+y_3)^2(x'_1+x'_3)\\
\alpha_3(x',y) & = & y_3-a(y_1+y_3)(x'_1+x'_3)\\
\alpha_4(x',y) & = & y_4+\left(cy_1-ay_2+(-b+2c)y_3+ay_4\right)x'_1-a(y_1+y_3)x'_2\\
& & \quad +\left((-b+2c)y_1-ay_2+(-2b+3c)y_3+ay_4\right)x'_3\\
& & \quad +a(y_1+y_3)x'_4+\ft(-b+c)(y_1+y_3)^2(x'_1+x'_3)\\
\beta_1(x',y) & = & x'_1-(x_1'+x'_3)(y_1+y_3)\\
\beta_2(x',y) & = & x'_2+(x_2'-x'_4)(y_1+y_3)+(x_1'+x'_3)(y_2-y_4)\\
& & \quad -\ft(-b+c)(x_1'+x'_3)^2(y_1+y_3)\\
\beta_3(x',y) & = & x'_3+(x_1'+x'_3)(y_1+y_3)\\
\beta_4(x',y) & = & x'_4+(x_2'-x'_4)(y_1+y_3)+(x_1'+x'_3)(y_2-y_4)\\
& & \quad -\ft(-b+c)(x_1'+x'_3)^2(y_1+y_3).
\end{eqnarray*}
We know from \S4 that $N_{\nabla,\nabla'}$ carries an invariant hypersymplectic structure whose associated neutral metric is complete (Theorem \ref{completef}). Moreover, using (\ref{R0}) and (\ref{ultima}), we may conclude that if $b=c$, then $\ngo_{\nabla,\nabla'}$ is 2-step nilpotent and the hypersymplectic metric is flat, whereas if $b\neq c$, then $\ngo_{\nabla,\nabla'}$ is 3-step nilpotent and the hypersymplectic metric is not flat. Note that taking $a,b,c\in\QQ$, the structure constants of $\ngo_{\nabla,\nabla'}$ with respect to the canonical basis $\{e_1,\ldots,e_4,f_1,\ldots f_4\}$ of $\ngo_{\nabla,\nabla'}$ are rational, and thus there exists a discrete co-compact subgroup of $N_{\nabla,\nabla'}$ \cite{M}. The complete non flat hypersymplectic metric on the Lie group induces a metric with the same properties on the associated compact quotient.

The Lie group $N_{\nabla,\nabla'}$ is diffeomorphic to $\RR^8$, hence there exists a global system of coordinates $x_1,\ldots,x_4,x'_1\ldots,x'_4$ such that the left-invariant 1-forms dual to the basis $\{e_1,\ldots,e_4,f_1,\ldots,f_4\}$ are given as follows
\begin{eqnarray*}
e^1 & = & \left(1-a(x'_1+x'_3)\right)\dif x_1-a(x'_1+x'_3)\dif x_3\\
e^2 & = & (-bx'_1+ax'_2-cx'_3-ax'_4)\dif x_1+\left(1+a(x'_1+x'_3)\right)\dif x_2\\
& & \quad +\left(-cx'_1+ax'_2-(-b+2c)x'_3-ax'_4\right)\dif x_3-a(x'_1+x'_3)\dif x_4\\
e^3 & = & a(x'_1+x'_3)\dif x_1+\left(1+a(x'_1+x'_3)\right)\dif x_3\\
e^4 & = & \left(-cx'_1+ax'_2-(-b+2c)x'_3-ax'_4\right)\dif x_1+a(x'_1+x'_3)\dif x_2\\
& &\quad +\left(-(-b+2c)x'_1+ax'_2-(-2b+3c)x'_3-ax'_4\right)\dif x_3+\left(1-a(x'_1+x'_3)\right)\dif x_4\\
f^1 & = & (x'_1+x'_3)\dif x_1+(x'_1+x'_3)\dif x_3+\dif x'_1\\
f^2 & = & \left(-\ft(-b+c)(x'_1+x'_3)^2-x'_2+x'_4\right)\dif x_1-(x'_1+x'_3)\dif x_2\\
& & \quad +\left(-\ft(-b+c)(x'_1+x'_3)^2-x'_2+x'_4\right)\dif x_3+(x'_1+x'_3)\dif x_4+\dif x'_2\\
f^3 & = & -(x'_1+x'_3)\dif x_1-(x'_1+x'_3)\dif x_3+\dif x'_3\\
f^4 & = & \left(-\ft(-b+c)(x'_1+x'_3)^2-x'_2+x'_4\right)\dif x_1-(x'_1+x'_3)\dif x_2\\
& & \quad +\left(-\ft(-b+c)(x'_1+x'_3)^2-x'_2+x'_4\right)\dif x_3+(x'_1+x'_3)\dif x_4+\dif x'_4\\
\end{eqnarray*}

The K\"ahler form $\omega_1$ on $\ngo_{\nabla,\nabla'}$ is $\omega_1=e^1\wedge e^2+e^3\wedge e^4+f^1\wedge f^2+f^3\wedge f^4$ and the hypersymplectic metric $g$ is given by
\[ g= -e^1\otimes f^2+e^2\otimes f^1-e^3\otimes f^4+e^4\otimes f^3+f^1\otimes e^2-f^2\otimes e^1+f^3\otimes e^4-f^4\otimes e^3.\]
In the particular case $a=0,\,b=1,\,c=0$, we obtain that
\begin{align*}
g = & \left(-\frac{3}{2}(x'_1+x'_3)^2+x'_2-x'_4\right)(\dif x_1+\dif x_3)^2 +2(x'_1+x'_3)(\dif x_1+\dif x_3)(\dif x_2-\dif x_4)-x'_1\dif x_1\dif x'_1\\
&  \quad -\dif x_1\dif x'_2+ \dif x_2\dif x'_1 + x'_3(\dif x_1\dif x'_3 +\dif x_3\dif x'_1)+ (x'_1 +2x'_3)\dif x_3\dif x'_3 -\dif x_3\dif x'_4+\dif x_4\dif x'_3
\end{align*}
is a complete non flat K\"ahler Ricci-flat neutral metric on $\RR^8$ .

\

\section{Final comments and questions}

We note that the complex structure $J$ defined in $\ngo_{\nabla,\nabla'}$ satisfies the condition $[Jx,Jy]=[x,y]$ for all $x,y\in \ngo_{\nabla,\nabla'}$, which implies the integrability of $J$. An almost complex structure $J$ on a Lie algebra $\ggo$ satisfying $[Jx,Jy]=[x,y]$ for all $x,y\in\ggo$ is called {\em abelian}. We also observe that $E$ satisfies a similar condition, $[Ex,Ey]=-[x,y]$ for all $x,y\in \ngo_{\nabla,\nabla'}$, which implies the integrability of $E$. An almost product structure $E$ on a Lie algebra $\ggo$ satisfying $[Ex,Ey]=-[x,y]$ for all $x,y\in\ggo$ will be called {\em abelian}. Related to these notions we have the following characterization:

\begin{prop}
Let $\{J,E\}$ be a complex product structure on the Lie algebra $\ggo$ and let $(\ggo,\ggo_+,\ggo_-)$ be the associated double Lie algebra, i.e., $\ggo_+$ and $\ggo_-$ are the Lie subalgebras of $\ggo$ such that $E|_{\ggo_+}=\Id,\,E|_{\ggo_-}=-\Id$. Then the following assertions are equivalent:
\begin{enumerate}
\item[\ri] $J$ is an abelian complex structure.
\item[\rii] The Lie subalgebras $\ggo_+$ and $\ggo_-$ are abelian;
\item[\riii] If $A_+$ and $A_-$ denote the annihilators of $\ggo_+$ and $\ggo_-$, respectively, in $\ggo^*$, then $\dif A_+\subset A_+\otimes A_-,\,\dif A_-\subset A_+\otimes A_-$, where $\dif:\ggo^*\longrightarrow\alt^2\ggo^*$ is given by $(\dif f)(x\wedge y)=-f([x,y])$.
\item[\riv] $E$ is an abelian product structure.
\end{enumerate}
\end{prop}

\begin{proof}
$\ri\Longleftrightarrow\rii$ Assume first that $J$ is abelian. If $x,y \in \ggo_+$ then $[x,y]\in \ggo_+$ and $[Jx,Jy]\in \ggo_-$ since $\ggo_+$ and $\ggo_-$ are subalgebras. But then $[x,y]=[Jx,Jy]\in \ggo_+ \cap \ggo_-=\{0\}$, and thus $[x,y]=[Jx,Jy]=0$ for all $x,y \in \ggo_+$. Thus, $\ggo_+$ and $\ggo_-$ are abelian.

Conversely, suppose that $\ggo_+$ and $\ggo_-$ are abelian. For $u,v\in \ggo_+$ or $u,v\in \ggo_-$, from the integrability of $J$ we obtain $[Ju,v]+[u,Jv]=J[u,v]-J[Ju,Jv]=0$. If $x=x_1+x_2,\, y=y_1+y_2$ with $x_1,y_1\in\ggo_+,\,x_2,y_2\in \ggo_-$, then
\begin{align*}
J([Jx,Jy]-[x,y]) & = [Jx_1+Jx_2,y_1+y_2]+[x_1+x_2,Jy_1+Jy_2] \\
                 & = [Jx_1,y_1]+[Jx_2,y_2]+[x_1,Jy_1]+[x_2,Jy_2] \\
                 & = ([Jx_1,y_1]+[x_1,Jy_1])+([Jx_2,y_2]+[x_2,Jy_2])\\
                 & = 0
\end{align*}
and thus $J$ is abelian.

\medskip

$\rii\Longleftrightarrow\riii$ Suppose first that $\ggo_+$ and $\ggo_-$ are abelian. Take $f\in A_+$, which is the annihilator of $\ggo_-$. It is known that $\dif f\in \alt^2A_+\oplus A_+\otimes A_-$ (see \cite{AS}), so we only have to see that the component of ($\dif f$) in $\alt^2A_+$ is zero. For $x,y\in\ggo_+$, we have
\[ (\dif f)(x\wedge y)=-f([x,y])=0, \]
showing that $\dif f\in A_+\otimes A_-$. The corresponding assertion for $f\in A_-$ follows in a similar manner.

Conversely, if $\riii$ is valid, take $f=f_1+f_2\in\ggo^*$ with $f_1\in A_+,\,f_2\in A_-$, and $x,y\in\ggo_+$ or $U,V\in\ggo_-$. Then
\[ f([x,y])=-(\dif f)(x\wedge y)= -(\dif f_1)(x\wedge y)-(\dif f_2)(x\wedge y)=0, \]
since $\dif f_1,\,\dif f_2\in A_+\otimes A_-$. Then $[x,y]=0$ and both $\ggo_+$ and $\ggo_-$ are abelian.

\medskip

$\rii\Longleftrightarrow\riv$ This follows by a straightforward computation.
\end{proof}

\medskip

If one of the conditions in the proposition above holds, we will say that the complex product structure $\{J,E\}$ is {\em abelian}. We will say that a hypersymplectic structure is abelian when the underlying complex product structure is abelian.

Using the previous proposition and Theorems 3.4 and 3.5 in \cite{AA} one can show that any Lie algebra carrying an abelian hypersymplectic structure is of the form given in Theorem 2.1., that is, a double product of two abelian Lie algebras endowed with compatible affine structures and symplectic forms. Furthermore, we showed in previous sections that in this case the Lie algebra is nilpotent and also the neutral metric is always complete and not necessarily flat.
It would be of interest to know geometric properties of neutral metrics compatible with $\{J,E\}$, $J$ and $E$ abelian, without imposing the condition on the associated forms being closed.

We also believe of interest to proceed as we did in this paper, carrying out the construction of double Lie groups from affine-symplectic data on another class of Lie groups (not necessarily abelian) and then understand the properties of the resulting hypersymplectic manifold.

\

\end{document}